\documentclass[12pt] {article}
\usepackage[dvips] {graphicx}
\usepackage{amssymb}
\usepackage{amsfonts}
\usepackage{amsmath}

\newtheorem{prop}{\bf Proposition}
\newtheorem{cor}{\bf Corollary}
\newtheorem{conj}{\bf Conjecture}
\newtheorem{thm}{\bf Theorem}

\def\half{\frac{1}{2}}

\def\ra{\rightarrow}
\def\reals{\mathbb{R}}
\def\ints{\mathbb{Z}}
\def\be{\mathbf{e}}
\def\csa{\mathfrak{h}}
\def\lag{\mathfrak{g}}
\def\g22{\Gamma^{2,2}}
\def\gr{\Gamma^{2,2}\otimes\mathbb{R}}

\begin{document}

\begin{titlepage}
\begin{flushright}
math.QA/0210451\\
DAMTP 2002-132\\[2cm]
\end{flushright}
\begin{center}{\Large \bf Lattice vertex algebras on general even,
self-dual lattices}\\[1.2cm]
Axel Kleinschmidt\\
Department of Applied Mathematics and Theoretical Physics\\
University of Cambridge, Wilberforce Road, Cambridge CB3 0WA, UK\\
Email: {\tt A.Kleinschmidt@damtp.cam.ac.uk}\\[0.5cm]

\end{center}

\renewcommand{\abstract}{\begin{center}\bf Abstract\\[.3cm]\end{center}}

\begin{abstract}
In this note we analyse the Lie algebras of physical states stemming
from lattice constructions on general even, self-dual lattices
$\Gamma^{p,q}$ with $p\ge q$.
It is
known that if the
lattice is at most Lorentzian, the resulting Lie
algebra is of generalized Kac-Moody type (or has a quotient that is). 
We show that this is not true as
soon as $q>1$. By studying a certain
sublattice in the case $q>1$ we obtain results that lead to the 
conclusion that the
resulting non-GKM Lie algebra cannot be described conveniently in terms of
generators and relations and belongs to a new and qualitatively
different class of Lie algebras.
\end{abstract}

\end{titlepage}
\begin{section}{Introduction}
\label{intro}

Vertex operators and their algebras have played an important r\^ole in
the development of string theory \cite{goddardolive1,frenkel} 
and constructions of algebras in
mathematics \cite{borcherds5,borcherds4,kac2}. The most prominent examples
have been chiral vertex algebras constructed from integral 
lattices which can
be interpreted as quantized (bosonic) strings moving in toroidal
spaces. Here the lattice is usually taken to be of Euclidean or Lorentzian
signature corresponding to compact Euclidean internal spaces or 
Minkowskian spacetimes respectively and the rank of the lattice in
this applicatio is
restricted by the no-ghost theorem \cite{goddardthorn}.

The formalisation and use of vertex algebra techniques eventually
led to the proof of the moonshine conjectures by Borcherds
\cite{borcherds4}. A major step on the road to this proof was the
introduction of the notion of a generalized Kac-Moody algebra (GKM for
short, see
\cite{borcherds2}). Algebras of this type  
can appear as Lie algebras derived from the
vertex operator construction and have several nice properties, most
notably the Weyl-Kac-Borcherds denominator formula

\begin{eqnarray} 
\label{den}
e^{\rho} \prod_{{\alpha}>0} (1-e^{\alpha})^{mult({\alpha})}
=\sum_{w\in W} \det (w)
w(e^{\rho}\sum_{\alpha}\epsilon({\alpha})e^{\alpha}),
\end{eqnarray}
which served to provide many new product formulae for known
automorphic forms \cite{borcherds1,borcherds3}.

Generalized Kac-Moody algebras\footnote{As the class of GKM algebras
comprises the usual finite-dimensional and Kac-Moody algebras we
restrict ourselves to discussing these.}
have also been proposed to underly several
structures in theoretical physics. They have been shown to control
threshold corrections in $N=2, d=4$ compactifications of the
heterotic string \cite{harveymoore1}, count degeneracies of certain
black holes \cite{dvv}, be relevant in Gromov-Witten theory
\cite{kawaiyoshioka} and possibly also as symmetries of M-theory
\cite{west1,juliapaulot}.

If the rank of the lattice is $\le 26$ and there is at most one
time-dimension, the no-ghost theorem from string theory implies that
the resulting Lie algebra is of GKM type or has a quotient that
is. For rank $>26$ arguments in \cite{borcherds6} show that the
algebra is still a GKM algebra.\footnote{I am grateful to
N. R. Scheithauer for bringing this point to my attention.}
Less is known about 
Lie algebras associated to lattices of signature $(p,q)$ where
both $p,q>1$. They have been proposed to be of relevance for the study
of algebras associated with BPS states in toroidal compactifications of
heterotic string theories \cite{harveymoore2}. 

This note aims at providing some results about the Lie algebras obtained
by lattice vertex operator constructions based on even, self-dual
lattices and to see how and where they differ from GKM type algebras. 
As the properties of the 
root systems of the algebras are intimately related to the
geometry of the lattices, we proceed by first studying systems of
fundamental roots of these lattices and then consider the implications
for the algebras. 

Our basic results deal with the existence of Weyl-like vectors and
suitable 
height functions for the lattice roots and we prove that neither
exist. 
This makes the existence of a nice set of fundamental roots
unlikely. We present some further evidence for this statement and
conjecture that the even, self-dual lattices of signature $(p,q)$ with
both $p,q>1$ do not possess a system of simple roots.

Similarly, the Lie algebra associated to the lattice at hand can then not be
written in terms of generators and relations in a fashion analogous to
the way one can construct GKM algebras. They thus must be part of a wider
class of Lie algebras with vastly different characteristics.

The structure of this note is as follows. To be mostly self-contained,
we first introduce some notation and terminology for the lattices
(section \ref{lattsec}), vertex algebras and GKM algebras 
(section \ref{valattsec}). We
indicate how the Lie algebra of physical states corresponding to the
lattice is constructed. In section \ref{reflrootsec} we discuss basic
properties of the lattice and algebra root systems and what we
understand by a set of fundamental or simple roots. The main results
mentioned above are contained in section \ref{qg1la}. Finally, we
offer some concluding remarks and possible further directions.

\end{section}

\begin{section}{Lattices}
\label{lattsec}

By an integral 
lattice $\Lambda$ of rank $n$ 
we understand the free Abelian group on a finite number
of generators $\be_i$ ($i=1,..,n$) together with a
non-degenerate, symmetric bilinear form
$(-,-):\,\Lambda\times\Lambda\ra\ints$ taking values in the
integers, so that

\begin{eqnarray}
\Lambda=\left\{ \sum_{i=1}^n n_i\be_i\, |\, n_i \in \ints\right\}.
\end{eqnarray}

Each such lattice can be seen as a discrete set of points in the
real vector space $\csa=\Lambda\otimes_\ints\reals$ which inherits an
inner product, also denoted $(-,-)$, from the bilinear form of the
lattice. This inner product is not necessarily definite. We can
associate to each integral lattice $\Lambda$ a dual lattice $\Lambda^*$
in $\csa$ by

\begin{eqnarray}
\Lambda^* :=\left\{\xi\in \csa\, | \, (\xi,\lambda)\in\ints
\;\mbox{for all }\;
\lambda\in\Lambda\right\}.
\end{eqnarray}

For integral lattices we obviously have $\Lambda\subset\Lambda^*$. 
We call a lattice {\sl even} if the norm squares of all its elements
are even, i.e. $\lambda^2=(\lambda,\lambda)\in 2\ints$ 
for all $\lambda\in\Lambda$. We call a
lattice {\sl self-dual} if it is identical to its dual lattice,
i.e. $\Lambda=\Lambda^*$. Self-duality for integral lattices is
equivalent to being {\sl uni-modular}, i.e. the inner product matrix
$g_{ij}=(\be_i,\be_j)$ 
has determinant of modulus one: $|\det (g_{ij})|=1$. 
We say $\Lambda$ has {\sl signature} $(p,q)$
if the inner product $(-,-)$ has $p$ positive and $q$ negative
eigenvalues. Necessarily we have $p+q=n$ and we assume $p\ge q$
without loss of generality.

We will restrict our attention in this paper to even, self-dual
lattices. It is a well-known result that such lattices exist if and 
only if
$(p-q)\equiv 0\mbox{mod }8$ and are unique if $q\ge 1$ \cite{serre}. These
lattices are usually denoted by $\Gamma^{p,q}$.

If $q=0$ we are
dealing with Euclidean even, self-dual lattices. There is one such of
rank $8$ (the root lattice $\Gamma^8$ of 
the Lie algebra $\mathfrak{e}_8$),
two such of rank $16$ ($\Gamma^8\oplus\Gamma^8$ and $\Gamma^{16}$,
the root lattice of  $\mathfrak{e}_8\times\mathfrak{e}_8$ and half the
weight lattice of
$\mathfrak{so}(32)$  
respectively). In $24$ dimensions there are $24$
even, self-dual lattices, the so-called Niemeier lattices among which the
Leech lattice is distinguished by having no elements of norm squared
equal to $2$.

If $q=1$ the lattices are said to be of Lorentzian signature and
sometimes also denoted by ${\rm II}_{n-1,1}$ with $n=2,10,18,26,\ldots$. The
most prominent of these are ${\rm II}_{1,1}$ (the root lattice of the
monster algebra \cite{borcherds4}), 
${\rm II}_{9,1}$ (the root lattice of the hyperbolic
Kac-Moody algebra $\mathfrak{e}_{10}$ \cite{gebnic}) and ${\rm II}_{25,1}$ (the root
lattice of the fake monster algebra \cite{borcherds1}).

If $q>1$ little seems to be known about these lattices in the
literature and it is part of the aim of this note to establish a few
basic properties of these lattices and their implications for the
vertex algebras based on them.

\end{section}

\begin{section}{Vertex Algebras of Lattices}
\label{valattsec}

We summarize the definition and the relevant properties of vertex
algebras and their construction starting from integral lattice,
including the associated Lie algebra of physical states. Relevant
references are \cite{goddardolive1,kac2,scheit1,frenkelaxiom}.

\begin{subsection}{Vertex Algebras}
\label{va}

A {\sl vertex algebra with conformal element} is a quadruple $({\cal
F},{\cal V},\omega,{\bf 1})$ consisting of a vector space $\cal F$,
called space of states, a map ${\cal V}: {\cal F}\rightarrow (End {\cal
F})[[z,z^{-1}]]$, and distinguished elements $\omega,{\bf 1}\in {\cal
F}$, termed the conformal vector and the vacuum respectively.
The map ${\cal V}$ is to be interpreted as associating to each
state $\psi\in {\cal F}$ a {\sl vertex operator} ${\cal
V}(\psi,z)=\sum_{n\in{\mathbb Z}}\psi_n z^{-n-1}$. One can think of a
state $\psi$ as being generated from the vacuum as $\psi=\lim_{z\rightarrow
0} {\cal V}(\psi,z){\bf 1}$. A vertex algebra
has to satisfy

\begin{enumerate}
\item{Regularity: $\psi_n\phi=0\;\;{\mbox{for $n$ sufficiently
large}}$}
\item{Vacuum: ${\bf 1}_n\psi=\delta_{n+1,0}\psi$}
\item{Injectivity: $\psi_n=0$ for all $n\in{\mathbb Z}\Leftrightarrow
\psi=0$}
\item{Conformal vector: $\omega_{n+1}=L_n$ and the $L_n$ obey the
Virasoro algebra and $L_0$ gives an integer grading of the space of
states
\begin{eqnarray}
{\cal F}=\oplus_{n\in {\mathbb Z}} {\cal F}_n\nonumber
\end{eqnarray}
with
\begin{eqnarray}
{\cal F}_n=\left\{\psi\in {\cal F}| L_0\psi=n\psi\right\}\nonumber
\end{eqnarray}
}
\item{Jacobi identity: 
\begin{eqnarray}
\sum_{i\ge 0}(-1)^i\left(\begin{array}{c}l\\i
\end{array}\right)\left(\psi_{l+m-i}(\phi_{n+i}\xi)- 
(-1)^l\phi_{l+n-i}(\psi_{m+i}\xi)\right)\nonumber\\
=\sum_{i\ge 0} 
\left(\begin{array}{c}m\\i\end{array}\right)
(\psi_{l+i}\phi)_{m+n-i}\xi
\end{eqnarray}}
\end{enumerate}

The Jacobi identity contains the most information about the
algebra. It comprises the ordinary Jacobi identity if
evaluated for $l=m=n=0$ and generalizations thereof.

A {\sl vertex operator algebra} is a vertex algebra for which 
the spectrum of
$L_0$ is bounded below and all ${\cal F}_n$ are
finite-dimensional. These conditions have obvious physical motivations. 
Interesting non-trivial examples are given by local
vertex operators in a conformal field theory as in the case of a
bosonic string compactified on an even and self-dual
lattice.

\end{subsection}

\begin{subsection}{Lattice construction}
\label{lattva}

For an integral lattice $\Lambda$ as in section \ref{lattsec} the
space $\csa/\Lambda$ is a torus and the vertex algebra we wish to
associate to the lattice $\Lambda$ is given by quantisation of the
space of maps from the circle $S^1$ into this torus. More precisely, we
consider the space

\begin{eqnarray}
\hat{\csa}=\csa[z,z^{-1}]\oplus\reals K
\end{eqnarray}
of Laurent polynomials in $V$ with central extension. 
We write $\alpha_{n,i}$ for $\be_i[z^n]$
and $\xi(n)=\xi^i\alpha_{n,i}$ for a general monomial $\xi[z^n]$ 
expanded in these basis vectors when $\xi=\xi^i\be_i$. Define a 
Heisenberg Lie algebra structure on this space by

\begin{eqnarray}
\left[\xi(m),\eta(n)\right]=m \delta_{m,-n}(\xi,\eta)K.
\end{eqnarray}

The space of states of the vertex algebra for the lattice $\Lambda$ is 

\begin{eqnarray}
{\cal F}_\Lambda= S(\csa^-)\otimes \reals\{\Lambda\},
\end{eqnarray}
where $\csa^-=z^{-1}\csa[z^{-1}]$ are the negative mode oscillators,
$S$ denotes the symmetric algebra and $\reals\{\Lambda\}$ is the
twisted group algebra of the lattice. 
Twisted here refers to the introduction of a
2-cocycle $\varepsilon:\;\Lambda\times\Lambda\ra\ints/2\ints$ 
which is necessary to obtain the right commutation relations, see
\cite{frenkel} for details. Elements
of the group algebra will be denoted by $e^\alpha$ with

\begin{eqnarray}
e^\alpha e^\beta=e^{\alpha+\beta}.
\end{eqnarray}

The Sugawara construction supplies the Virasoro element by

\begin{eqnarray}
\omega=\half\sum_{i=1}^n h^i(-1)h_i(-1)\otimes e^0,
\end{eqnarray}
where $h^i,h_i$ are dual bases of $\csa$. The vacuum is given by
$1\otimes e^0$.

The map $\cal V$ associating a vertex operator to each state is
defined in the usual way using the normal ordering of the Heisenberg
modes. It can be shown that with these choices we have constructed a
vertex algebra with conformal element \cite{kac2}. 

\end{subsection}

\begin{subsection}{Lie algebra of physical states}
\label{laphys}

We can use the conformal structure of the lattice vertex algebra 
to define a Lie algebra structure
on the space of physical states. Define

\begin{eqnarray}
P_n:=\left\{\psi\in {\cal F}_n: L_m\psi=0\;\mbox{for }\;m>0\right\}
\end{eqnarray}
to be the space of physical states with weight $n$. Then it can be
shown \cite{borcherds5} that the space $P_1/L_{-1}P_0$
carries the structure of a Lie algebra with respect to the product

\begin{eqnarray}
[\psi,\phi]=\psi_0\phi=Res_{z=0}{\cal V}(\psi,z)\phi.
\end{eqnarray}

This construction works for all vertex algebras with conformal vector
but we will only consider the case of vertex algebras derived from
a lattice construction for even, self-dual lattices
$\Gamma^{p,q}$. 
In this case
the root lattice of the Lie algebra of physical states is just
$\Gamma^{p,q}$ and the roots are the non-zero 
elements with norm squared less than
$2$ and the multiplicity of a root $\alpha\in\Gamma^{p,q}$ is given by
the classical partition function

\begin{eqnarray}
\mbox{mult}(\alpha)=p_{d-1}(1-\half\alpha^2)-p_{d-1}(-\half\alpha^2),
\end{eqnarray}
where the two terms describe the number of independent elements for
given norm squared in $P_1$ and $P_0$. $p_k(M)$ is the number of
partitions of an integer $M$ into positive integers with $k$ colours
and is generated by

\begin{eqnarray}
\sum_{m\ge 0} p_k(m)q^m=\prod_{m\ge 1}(1-q^m)^{-k}.
\end{eqnarray}
The Cartan subalgebra is given by all states with zero momentum and
arbitrary polarization vector in $\Gamma^{p,q}\otimes\reals$.

If this construction is applied to the lattice ${\rm II}_{25,1}$ one arrives
at the fake monster algebra \cite{borcherds1} which is an example of
a generalized Kac-Moody algebra (which will be defined in the
following section \ref{gkm}) and in general for a Lorentzian lattice
we obtain a GKM algebra \cite{borcherds6}. 
So the Lie
algebras for $\Gamma^{p,q}$ with $q\le 1$ and are at most
of generalized Kac-Moody type, often they are even
finite-dimensional or of Kac-Moody type.

\end{subsection}

\begin{subsection}{Generalized Kac-Moody Algebras}
\label{gkm}

Following \cite{borcherds2}, we briefly review a few standard facts
about generalized Kac-Moody algebras. A {\sl
generalized Kac-Moody algebra} is defined via a symmetrized Cartan
matrix $A=(a_{ij})$, with $i,j$ in a possibly infinite but countable
index set $I$, satisfying the following conditions

\begin{eqnarray}
\label{gcmprops}
a_{ij}=a_{ji},\nonumber\\
a_{ij}\le 0 \;\;{\mbox{if $i\ne j$}},\\
\frac{2 a_{ij}}{a_{ii}}\in {\mathbb{Z}}\;\; {\mbox{if $a_{ii}>
0$}}.\nonumber
\end{eqnarray}

We also assume the existence of a real vector space $H$ with symmetric
bilinear inner product (not necessarily positive definite) and
elements $h_i\in H$ ($i\in I$) such that $(h_i,h_j)=a_{ij}$. Then we
define the generalized Kac-Moody algebra $G$ to be the Lie algebra
generated by H and $e_i$ and $f_i$ subject to the relations

\begin{enumerate}
\item{The image of $H$ in $G$ is commutative.}
\item{$h\in H$ acts diagonally on the $e_i,f_i$: $[h,e_i]=(h,h_i)e_i$
and $[h,f_i]=-(h,h_i)f_i$.}
\item{$[e_i,f_i]=h_i$, $[e_i,f_j]=0$ if $i\ne j$.}
\item{If $a_{ii}>0$ then $(\text{ad}\, e_i)^{1-2a_{ij}/a_{ii}}e_j=0$ and
$(\text{ad}\, f_i)^{1-2a_{ij}/a_{ii}}f_j=0$ (Serre relations).}
\item{If $a_{ij}=0$ then $[e_i,e_j]=[f_i,f_j]=0$.}
\end{enumerate}

The main difference compared to finite-dimensional simple Lie algebras or
Kac-Moody algebras is that imaginary simple roots are permitted as
$a_{ii}$ can be less or equal to $0$. The multiplicity of an imaginary
simple root can be greater than $1$. Also, the number of simple
roots can be infinite dimensional. Most terminology and properties
carry over. We extend the bilinear form on $H$ to the root lattice
$Q$ which is the free Abelian group on the simple elements. 
The Weyl group is generated by reflections in the real 
simple roots, and roots can be grouped into positive or negative ones
according to whether they are a sum of simple roots or the negative
thereof.

Among the standard properties of GKM algebras we mention the
following. The denominator formula (\ref{den}) involving the Weyl
vector determined by the real simple roots holds. Unless for
extremely degenerate cases, a GKM algebra $G$
has an integer grading $G=\oplus G_m$ 
such that each $G_m$ is finite-dimensional
if $m\ne 0$ and $G_0$ is Abelian.\footnote{This is true unless the
Cartan matrix has a infinite number of identical rows.} Similar to Kac-Moody algebras, 
$G$ has an invariant bilinear form such that $G_\alpha$ and
$G_\beta$ are orthogonal unless $\alpha+\beta=0$, and the Hermitian form
constructed from this by use of the Chevalley involution is positive
definite on the root spaces.

\end{subsection}

\end{section}

\begin{section}{Reflections, Algebra and Lattice Roots}
\label{reflrootsec}

In order to analyse the Lie algebra of physical states we construct
from lattices further we need to understand the root system and the
geometry of the root lattice. In this section we introduce the
required notation.

\begin{subsection}{Lattice Roots and Algebra Roots}
\label{lattroots}

Following \cite{conwaysloane}, we call an 
element $\alpha\in\Gamma^{p,q}$ a {\sl lattice root} if the
associated elementary reflection 

\begin{eqnarray}
w_\alpha(\gamma)=\gamma-2
\frac{(\gamma,\alpha)}{(\alpha,\alpha)}\alpha,\;\;\gamma\in\Gamma^{p,q}
\end{eqnarray}
is a symmetry of the lattice. It is not difficult to see that for
even, self-dual lattices the lattice roots are precisely the elements
of norm squared $2$: 

\begin{eqnarray}
\Delta^{latt}=\left\{\alpha\in\Gamma^{p,q}\;|\; \alpha^2=2\right\}.
\end{eqnarray}
The set of all reflection forms a group $\cal R$, called the reflection
subgroup of the automorphism group of the lattice.\\

If $\mathfrak{g}$ is a Lie algebra with root lattice $\Gamma^{p,q}$
and Cartan subalgebra $\csa$ then it possesses a natural grading

\begin{eqnarray}
\mathfrak{g}=\bigoplus_{\alpha\in\Gamma^{p,q}}\mathfrak{g}_\alpha
=\csa\oplus\bigoplus_{0\ne\alpha\in\Gamma^{p,q}}\lag_\alpha,
\end{eqnarray}
where

\begin{eqnarray}
\mathfrak{g}_\alpha=\left\{x\in\mathfrak{g}\;|\;
[h,x]=\alpha(h)x\;\mbox{for all }\,h\in\csa\right\}.
\end{eqnarray}
If $\mathfrak{g}_\alpha\ne\left\{0\right\}$ and $\alpha\ne 0$ 
then $\alpha$ is called an {\sl algebra root} or root for short and
$\mathfrak{g}_\alpha$ the associated root space. If $\alpha^2>0$,
$\alpha$ is called {\sl real} and {\sl imaginary} otherwise, denoted
by $\Delta^{re}$ and $\Delta^{im}$ respectively. For the
lattice Lie algebras we are interested in all real roots have norm
squared $2$. In fact, here 
all $\alpha$ with $\alpha^2=2$ are algebra roots
and 

\begin{eqnarray}
\Delta^{latt}=\Delta^{re}
\end{eqnarray}
for these lattice Lie algebras. This is not true in general, for
instance for most Kac-Moody algebras. In a similar fashion to above,
the real roots have reflection symmetries which form a group $\cal W$
called the Weyl group of the algebra. We note that the set of all
roots is invariant under $\alpha\ra-\alpha$, also reflected in
the
Chevalley involution of the Lie algebra.

\end{subsection}

\newpage

\begin{subsection}{Fundamental and simple roots}
\label{fundsim}

It is an important question if the set of all roots can be
conveniently described
in terms of some generators as for GKM type algebras. As upon commutation
$[\lag_\alpha,\lag_\beta]\subset\lag_{\alpha+\beta}$ an element which
can be written as the sum (or commutator) of other elements cannot be
part of a minimal generating set. It is natural to
look for such a generating set by a minimality condition with respect to
addition. For this it is first necessary to define the notion of
positivity and negativity of a (lattice or algebra) root which
respects the group operation on the lattice. This can be done, for
instance, 
by finding a hyperplane in the vector space
underlying the lattice not containing any of the roots (real or
imaginary). Then call one side of the hyperplane positive, the other
one negative. 

For the case of the lattice roots, those positive roots which cannot
be written as a sum of other positive roots will be called {\sl 
fundamental
roots} and are denoted $(\alpha_i)$ for $i$ in some index 
set $I$.\footnote{This notion is not fully equivalent to saying that the
fundamental roots are orthogonal to the faces of a fundamental domain
for the action of ${\cal R}$ on $\csa$. The image of the fundamental
chamber $\cap_{i\in I} \left\{\xi\in\csa\;|\;(\alpha_i,\xi)\ge
0\right\}$ under $\cal R$ 
is the so-called Tits cone $X$. An example when $X\ne\csa$
is given by ${\rm II}_{9,1}$. Nevertheless, $X$ contains all lattice roots.}
The set of fundamental roots is sufficient to describe the reflection
group ${\cal R}$ and the information can be encoded in a Dynkin
diagram or a Cartan-like matrix.

For the lattice Lie algebra, those positive roots whose root spaces
contain elements which cannot be
written as a commutator of elements of the root spaces of 
other positive roots will be called {\sl simple roots} with multiplicity
corresponding to the number of linearly independent such elements in their root
spaces. The fundamental roots of the lattice will turn into simple
real roots of the algebra and then in principle one might try to
determine the additional simple roots recursively if one has a
suitable grading like the distance to the hyperplane chosen for
splitting the roots into positive and negative. Again, the simple
roots will be denoted $\alpha_i$ ($i\in I$) and their inner-product
matrix contains all information about the Lie algebra. 

It seems an open question under which circumstances such a generating
set exists and can actually be determined.

\end{subsection}

\begin{subsection}{Weyl-like elements}
\label{weylvec}

An important quantity for generalized Kac-Moody algebras is the Weyl
element $\rho$ which is defined by 

\begin{eqnarray}
(\rho,\alpha_i)=1
\end{eqnarray}
for all real simple roots $\alpha_i$ where one might need to extend
$\csa$ by derivations if the Cartan matrix is singular. It 
derives its importance from the r\^ole it plays in the character
formula and, in light of the discussion in the preceding paragraph,
can be seen to be an optimal choice of grading for the positivity
problem\footnote{For this we might require that the lattice roots 
span the lattice.}.
This said, we call an element $\rho\in\Gamma^{p,q}$ {\sl Weyl-like} if
it has non-vanishing inner product with all lattice roots:

\begin{eqnarray}
(\rho,\alpha)\ne 0\;,\;\mbox{for all}\;\alpha\in\Delta^{latt}.
\end{eqnarray}
Such an element can then be used to find the fundamental roots of the
lattice.\footnote{For most purposes one might also require that to
each integer $m$ there are only finitely many roots which have inner
product $m$ with $\rho$. That this is not strictly necessary can be
seen from the example of ${\rm II}_{25,1}$.}

\end{subsection}

\end{section}

\begin{section}{Properties of the Lie algebra of physical states based
on the lattice $\Gamma^{p,q}$ for $q>1$}
\label{qg1la}

Whereas the Lie algebra $\lag$ of physical states for at
most Lorentzian lattices can be classified as
generalized Kac-Moody algebras (GKM), we will see that the 
situation is different as soon as
$q>1$.

The lattice $\Gamma^{p,q}$ can be written as

\begin{eqnarray}
\Gamma^{p,q}=\Gamma^{8s}\oplus 
\underbrace{{\rm II}_{1,1}\oplus\ldots\oplus{\rm II}_{1,1}}_q,
\end{eqnarray}
where $\Gamma^{8s}$ is a Euclidean lattice of rank $8s$ where
$s=(p-q)/8$ and there are $q$ summands of the basic light-cone lattice
${\rm II}_{1,1}$. For $q>1$, $\Gamma^{p,q}$ contains at least one copy
of 
$\Gamma^{2,2}={\rm II}_{1,1}\oplus{\rm II}_{1,1}$ as a sublattice. Many of the
problems in analysing $\lag$ based on $\Gamma^{p,q}$ can already be
understood by studying $\Gamma^{2,2}$.

A convenient description of the lattice $\g22$ is given 
in terms of a double set of light-cone co-ordinates:

\begin{eqnarray}
\Gamma^{2,2}=\left\{
(k,l;m,n)\;|\;k,l,m,n\in\mathbb{Z}\right\},
\end{eqnarray}
equipped with the metric

\begin{eqnarray}
|(k,l;m,n)|^2=-2kl-2mn.
\end{eqnarray}

As discussed above. the lattice roots 
are exactly those elements of $\g22$ with $\alpha^2=2$ or,
put differently, the integral solutions to the non-linear 
diophantine problem

\begin{eqnarray}
\label{rootcond}
kl+mn=-1.
\end{eqnarray}
There is an obvious one-to-one correspondence between the lattice
roots and elements of $M\in SL(2,\ints)$ by making use of the
bijection  $k\leftrightarrow -k$ and multiplying (\ref{rootcond}) by $-1$.

Now we are able to show some interesting properties of the lattice.

\begin{prop}
\label{weyllike}
There is no Weyl-like element for the lattice $\Gamma^{2,2}$.
\end{prop}

{\bf Proof:} Assuming that there is such an element $\rho=(a,b;c,d)$ which can
be written as a $2\times 2$-matrix with integer entries

\begin{eqnarray}
A=\left(\begin{array}{cc}-b&c\\d&a\end{array}\right),
\end{eqnarray}
we have to
show that there always exists an element $R\in SL(2,\ints)$ 

\begin{eqnarray}
R=\left(\begin{array}{cc}-k&m\\n&l\end{array}\right)
\end{eqnarray}
(corresponding to a lattice root $\alpha=(k,l;m,n)$) such that

\begin{eqnarray}
\mbox{tr}(RA)=al+bk+cn+dm=-(\rho,\alpha)=0.
\end{eqnarray}
The condition $\mbox{tr}(RA)=0$ for some $R\in SL(2,\ints)$ is
equivalent to the condition

\begin{eqnarray}
\mbox{tr}(PAQ)=0,
\end{eqnarray}
for some matrices $P,Q\in SL(2,\ints)$. Hence, if we can bring $A$ to a
traceless form by elementary row and column manipulations we are done. But by a
variant of Euclid's algorithm we can always arrive at the following form for 
$A$ by doing such manipulations

\begin{eqnarray}
\left(\begin{array}{cc}d_1&0\\0&d_1e\end{array}\right),
\end{eqnarray}
for some $d_1\ne 0$ (unless $A=0$ to start with) and
$e\in\ints$. Now, we can add the first column to the second and then
subtract $e+1$ times the first row off the second and the resulting
matrix will be traceless.$\square$\\

The significance of proposition \ref{weyllike} is that we cannot find
a linear ``height'' function $h:\;\g22\ra\ints$ which is non-vanishing on the
lattice 
roots and can be used to distinguish positive and negative roots and
to subsequently determine the fundamental roots. Conversely, supposing
we have somehow arrived at a set of fundamental roots, so that all
lattice roots can be written either as a sum
or minus a sum of fundamental roots, there cannot be a Weyl vector for
this set of fundamental roots in $\g22$. Proposition \ref{weyllike} 
shows that $\g22$ cannot be graded by integers in a fashion desirable
for GKM algebras.

\begin{prop}
\label{genht}
Every group homomorphism $h:\;\g22\ra\reals$ which is non-vanishing on
$\Delta^{latt}$ has zero as an accumulation point of its image.
\end{prop}

{\bf Proof:} By proposition \ref{weyllike} we only need to consider
the case when $h$ takes irrational values, at least for some lattice
roots. Let $(k_0,l;m_0,n)$ be such an element. We know that $l,n$ are
coprime and that for all $t\in\ints$ 

\begin{eqnarray}
(k_0+t n,l;m_0-t l,n)
\end{eqnarray}
is also a root\footnote{Actually, all solutions to the lattice root
condition are of this form: For every coprime pair $l,n$ there is such
a one-parameter family of solutions.}. If $h$ is given by 

\begin{eqnarray}
h(k,l;m,n)=-al-bk-cn-dn,
\end{eqnarray}
let

\begin{eqnarray}
h_{l,n}(t)&=&-al-b(k_0+tn)-cn-d(m_0-tl)=(dl-bn)t-al-bk_0-cn-dm_0\nonumber\\
&=&Mt+N.
\end{eqnarray}
There exists a $t_0\in\ints$ such that $h_{l,n}(t_0)\in(0,M)$. Now
approximate $b$ and $d$ over the rationals as

\begin{eqnarray}
b\approx\frac{p_1}{q_1},\;d\approx\frac{p_2}{q_2}
\end{eqnarray}
with $q_1p_2$ and $q_2p_1$ coprime. Choosing $l=q_2p_1$ and $n=q_1p_2$
gives $M=dl-bn\approx 0$. By increasing the accuracy of the
approximation we can make $M$ arbitrarily small, showing that $0$ is
an accumulation point of $h(\Delta^{latt})$.$\square$\\

Thus, even though it is possible to split the lattice roots into
positive and negative roots by a general hyperplane there is never a
minimal element that could be identified as fundamental. Combining
proposition \ref{weyllike} and proposition \ref{genht}, we 
get

\begin{cor}
\label{hyperplanes}
Let $H\subset\g22\otimes\reals$ be a hyperplane. Then either $H$
contains an element of $\Delta^{latt}$ or there is an element of
$\Delta^{latt}$ arbitrarily close nearby.
\end{cor}

Interpreting the possible maps $h$ as points in $\g22\otimes\reals$ we
also get that to each point in $\g22\otimes\reals$ there is a
hyperplane orthogonal to a lattice root passing arbitrarily close
by. This has implications for the reflection group $\cal R$.

\begin{cor}
\label{funddom}
Any potential fundamental domain of the 
group $\cal R$ generated by reflections in the lattice roots of
$\g22$ acting on $\g22\otimes\reals$ has empty
interior. 
\end{cor}

{\bf Proof:} Splitting
$\Delta^{latt}=\Delta^{latt}_+\cup\Delta^{latt}_-$ every single
reflection in an element $\alpha\in\Delta^{latt}_+$ exchanges the two
half spaces

\begin{eqnarray}
H_\alpha^+&=&\left\{\gamma\in\gr |\;(\alpha,\gamma)\ge
0\right\}\;\;\mbox{and}\nonumber\\
H_\alpha^-&=&\left\{\gamma\in\gr |\;(\alpha,\gamma)\le 0\right\},
\end{eqnarray}
leaving the orthogonal elements invariant. We pick the positive half
space and try to construct a fundamental domain $D$ for the action of $\cal R$
on $\gr$.

\begin{eqnarray}
D=\bigcap_{\alpha\in\Delta^{latt}_+}H_\alpha^+
\end{eqnarray}
is a closed set with the property that every reflection takes a point
in the interior of $D$ to a point outside of $D$ and $D$ is bounded by
hyperplanes orthogonal to the lattice roots. By corollary
\ref{hyperplanes}, any point in $\gr$ has such a hyperplane nearby
and thus the
interior of $D$ has to be empty.$\square$\\

This also relates to the absence of a Weyl vector which would lie in
the interior of such a set and a ``ball'' of radius $1$ around its tip
touching all the faces. We can describe $D$ more explicitly when we
use a height function $h$ as above. Then $D$ is a ray
in $\gr$, precisely the one in the direction of the point determining
$h$ and the Tits cone is $(\gr)\backslash\g22$.

We have accumulated several pieces of evidence for

\begin{conj}
\label{nofundroots}
The lattice $\g22$ does not possess a system of fundamental roots.
\end{conj}

We remark that choosing non-linear orderings of the lattice roots like
lexicographic ordering one can find a few fundamental elements which
can be shown not be a generating set of all roots. Completing them to
a set of fundamental roots cannot be done algorithmically as there is
no analogue of height, i.e. a superlinear functional that ensures that
one only has to consider a given set of roots when asking if a
particular root can be written as a sum of ``lower'' roots. 
As the fundamental roots of the lattice carry over by construction to
the real simple roots of the Lie algebra of physical states $\lag$. As
$\g22\subset\Gamma^{p,q}$ for $p\ge q>1$ we are led to the believe
that $\lag$ cannot have a system of simple roots analogous to the
construction of generalized Kac-Moody algebras.

Using the results on the lattice $\g22$ we can now prove

\begin{prop} 
\label{notgkm}
The Lie algebra of physical states $\lag=P_1/L_{(-1)}P_0$
associated to the lattice $\Gamma^{p,q}$ is not of generalized
Kac-Moody type for $q>1$.
\end{prop}

{\bf Proof:} We restrict to the subalgebra where we only consider the
lattice $\g22$ and assume that $\lag$ is 
a GKM algebra and show that this leads to a
contradiction. So suppose $\lag$ is generated by a set of simple roots
$\{\alpha_i:i\in I\subset\mathbb{N}\}=\Pi=\Pi^{re}\uplus\Pi^{im}$ 
with generalized Cartan matrix
$a_{ij}=(\alpha_i,\alpha_j)$ and we can exclude zero rows. 
We denote the set of roots of $\lag$ as $\Delta^Q\subset Q$ which
splits into positive ($\alpha>0$) and negative ($\alpha<0$) roots
$\Delta^Q=\Delta^Q_+\uplus\Delta^Q_-\subset Q_+\uplus Q_-$. We define the fundamental set
$C=\{\alpha\in Q_+| (\alpha,\alpha_i)\le 0$ for all real simple 
$\alpha_i$ and the support of $\alpha$ is connected$\} \backslash
\bigcup_{j\ge 2}j\Pi^{im}$.

We know that our algebra is graded by the lattice $\g22$ as
$\lag=\bigoplus_{\alpha\in\g22}\lag_\alpha$,
where each $\lag_\alpha$ is finite dimensional and we are abusing
notation by denoting roots in $Q$ and $\g22$ by $\alpha$.
If $\lag$ is a GKM algebra then $\gr$ will be the quotient of
$Q\otimes\mathbb{R}$ by the
kernel $K$ of the Cartan matrix. 
Suppose the kernel contains a (positive) 
root $\alpha$ of $\lag$
then the norm of this root is $\alpha^2=0$ and thus 
$\alpha$ is an
isotropic root in the Weyl chamber. 
As we excluded zero rows the support of $\alpha$ is
affine. For a simple root $\beta$ from the support of $\alpha$, the
elements $\beta+ n \alpha$ would be roots for all natural $n$, giving
rise to infinite multiplicities in the quotient contradicting what we
know about $\lag$. Hence $K$ does not contain any root.

Now suppose there are elements $\alpha>0$ and $\beta<0$ which differ
by an element $\kappa$ of the kernel $K$,
i.e. $\alpha=\beta+\kappa$. This $\kappa$ obviously has to be positive
and lies in the Weyl chamber.
As it is not a root it has to have
disconnected support in order not to lie in the fundamental set.
So it is the sum of null roots of disconnected
subdiagrams of the Dynkin diagram of $\lag$. But then each of the
summands is also in the kernel (by disconnectedness) and thus a null
root of $\lag$ which contradicts the fact above that there are no roots
in the kernel. Hence $K\cap Q_\pm=K\cap\Delta^Q_\pm=\{0\}$. This means
that the quotient map takes $Q_\pm$ to disjoint, convex sets $\g22_\pm$ which
by linearity are exchanged by $p\leftrightarrow -p$ and thus can be
split by a hyperplane in $\g22$. By the properties of $\g22$ such a
hyperplane must have a normal vector $\rho$ with irrational coordinates and
the fundamental domain is the ray in the direction of $\rho$.
This is where all imaginary simple roots must live but as the
coordinates of $\rho$ are irrational this ray does not contain any
lattice points and so there are no imaginary simple roots for this Lie
algebra. So $\lag$ can be at
most a standard Kac-Moody algebra \cite{kac1} with an infinite number
of (real) simple roots.

The set $\Delta^{im}_+$ ($\Delta^{im}_-$) is Weyl invariant so that any sum of
positive (negative) imaginary roots
roots will again be an imaginary root if it is a root. If we
separate the positive and negative roots in
$\g22$ by $\rho=(a,b;c,d)$, we can assume that $b\ne 0$ without loss
of generality as we can swap $a$ and $b$ by a Weyl reflection and not
both can vanish. Now consider the roots
$\beta_1=(N,0;1,0)$ and $\beta_2=(N,0;0,-1)$ which are both imaginary
and if we choose $N$ large enough both are either positive or
negative. Their sum is then also positive or negative but a real root
contradicting the structure on $\Delta^{im}_\pm$ explained above.$\square$\\

Hence, $\lag$ is an element of a wider class of Lie algebras. 
As the defining relations of generalized Kac-Moody algebras are
automatically satisfied due to the vertex operator construction, we
observe that if there were simple roots for $\lag$ they they would
need to violate at least one of the basic properties of generalized
Cartan matrices (\ref{gcmprops}). 
Closer inspection reveals that their inner product
matrix would need to have positive off-diagonal entries. For such
Lie algebras we can work out the multiplicities for any given root
inductively by height using the theory of generalized Kac-Moody
algebras. 

\begin{thm}
\label{nosimroots}
For $q>1$, the Lie algebra $\lag$ associated to $\Gamma^{p,q}$ 
cannot be described in terms of simple generators and relations
analogous to GKM algebras, say, encoded in Cartan matrices or
Dynkin-like diagrams.
\end{thm}

{\bf Proof:} Assume there is a system of generators for $\lag$ and as
argued above they satisfy by construction relations similar to the
ones of GKM algebras. But their Cartan matrix cannot be that of a GKM
algebra as we have shown that this leads to a contradiction. So we
have to relax the property that we only allow non-positive
off-diagonal elements in the Cartan matrix. But for the elements used
above we can view these as elements of an embedded 
GKM algebra (if we restrict
the height accordingly) where we can
still obtain the same contradiction as in the proof of proposition
(\ref{notgkm}). $\square$\\

\end{section}

\begin{section}{Conclusions and Open Questions}
\label{concl}

The limitations on the type of Lie algebra one can obtain for this
general lattice construction seem rather tight. Is there any other way to
describe these algebras conveniently except for saying that they
derive from the lattice construction? For instance, the class of
algebras excluded here rest on amalgamating $\mathfrak{sl}_2$
subalgebras -- maybe if one took other basic objects one can construct
these Lie algebras. Is there some generalization of
a character or denominator formula one could write down? As there is
no Weyl-like vector and also the Weyl group seems rather unwieldy, it
is not clear how to suitably replace the quantities appearing in the
denominator formula (\ref{den}). The representation theory for these
algebras should also prove interesting.

Besides these mathematical questions it is not clear what the
relevance of these algebras could be in theoretical physics, along the
lines of \cite{harveymoore1,harveymoore2}. If algebras of BPS states
turn out to be similar to these general lattice Lie algebras it is
essential to get some handle on the underlying
structure. \\

{\bf Acknowledgements}\newline
I would like to thank my supervisor P. Goddard and M. R. Gaberdiel for
valuable support and stimulating advice. I am also grateful to
T. Fisher, I. Grojnowski, D. Olive and 
S. Sch\"afer-Nameki for discussions and to N. R. Scheithauer for
comments. Support by
the Studienstiftung des deutschen Volkes, EPSRC,
and the Cambridge European Trust is gratefully acknowledged.

\end{section}


\begin{thebibliography}{30}

\bibitem{goddardolive1} P. Goddard and D. Olive, {\sl
Algebras, Lattices, and Strings} in ``Vertex Operators in mathematics
and physics'', ed. J. Lepowsky et al., Springer Verlag 1985

\bibitem{frenkel} I. B. Frenkel, {\sl Representations of Kac-Moody
algebras and dual resonance models} in ``Applications of Group Theory
in Physics and Mathematical Physics'', Vol. 21, Lectures in Appplied
Mathematics, M. Flato, P. Sally, G. Zuckerman, eds. AMS 1985

\bibitem{borcherds5} R. E. Borcherds, {\sl Vertex algebras, Kac-Moody
algebras and the monster}, Proc. Nat. Acad. Sci. U.S.A. {\bf 83}
(1986), 3068-3071

\bibitem{borcherds4} R. E. Borcherds, {\sl Monstrous moonshine and
monstrous Lie superalgebras}, Invent. Math.  {\bf 109} (1992), 405-444

\bibitem{kac2} V. Kac, {\sl Vertex Algebras for Beginners}, 2nd
edition, AMS, Providence, Rhode Island, 1998

\bibitem{goddardthorn} P. Goddard and C. B. Thorn, {\sl Compatibility
of the dual Pomeron with unitarity and the absence of ghosts in the
dual resonance model}, Phys. Lett. {\bf 43} No. 2 (1972), 235-238

\bibitem{borcherds2} R. E. Borcherds, {\sl Generalized Kac-Moody
Algebras}, J. Algebra {\bf 115} (1988) 501-512

\bibitem{borcherds1} R. E. Borcherds,  {\sl The
monster Lie algebra}, Adv. Math. {\bf 83}, No. 1 (1990) 30-47

\bibitem{borcherds3} R. E. Borcherds, {\sl Automorphic Forms on
$O_{s+2,2}({\bf R})$ and infinite products}, Invent. Math. {\bf 120}
(1995) 161-213

\bibitem{harveymoore1} J. A. Harvey and G. Moore, {\sl Algebras, BPS
states, and Strings}, Nucl. Phys. {\bf B 463} (1996) 315-368,
{\tt hep-th/9510182}

\bibitem{dvv} R. Dijkgraaf, E. Verlinde and H. Verlinde, {\sl
Counting Dyons in $N=4$ String Theory}, Nucl. Phys {\bf B 484} (1997)
543-561, {\tt hep-th/9607026}

\bibitem{kawaiyoshioka} T. Kawai and K. Yoshioka, {\sl String Partition
Functions and Infinite Products}, Adv. Theor. Math. Phys. {\bf 4}
(2001) 397-485, {\tt hep-th/0002169 }

\bibitem{west1} P. West, {\sl $E_{11}$ and M Theory},
Class. Quant. Grav. {\bf 18} (2001) 4443-4460, {\tt hep-th/0104081}

\bibitem{juliapaulot} P. Henry-Labord\`ere, B. Julia, L. Paulot, {\sl
Borcherds symmetries in M-theory}, JHEP 0204 (2002) 049, {\tt
hep-th/0203070}

\bibitem{borcherds6} R. E. Borcherds, {\sl A characterization of
generalized Kac-Moody algebras}, J. Algebra {\bf 174} (1995) 1073-1079

\bibitem{harveymoore2} J. A. Harvey and G. Moore, {\sl On
the algebras of BPS states}, Commun. Math. Phys. {\bf 197} (1998)
489-519, {\tt hep-th/9609017}

\bibitem{serre}  J.-P. Serre, {\sl A Course in Arithmetic},
Springer-Verlag, New York, 1973

\bibitem{gebnic} R. W. Gebert and H. Nicolai, {\sl On $E_{10}$ and the
DDF Construction}, Commun. Math. Phys. {\bf 172} (1995) 571-622,
hep-th/9406175

\bibitem{scheit1} N. R. Scheithauer, {\sl Vertex algebras, Lie
algebras, and superstrings}, J. Algebra {\bf 200} No. 2 (1998), 
363-402, {\tt hep-th/9802058}

\bibitem{frenkelaxiom} I. B. Frenkel, Y. Huang and J. Lepowsky, {\sl
On axiomatix approaches to vertex operator algebras and modules},
Mem. Amer. Math. Soc. {\bf 104} (1993)

\bibitem{conwaysloane} J. H. Conway and N. J. A. Sloane, {\sl Sphere
Packings, Lattices and Groups}, Springer-Verlag, New York, 1988

\bibitem{kac1} V. G. Kac, {\sl Infinite dimensional Lie algebras}, 3rd
edition, Cambridge University Press (1990)

\end{thebibliography}
\end{document}